\documentclass[a4paper,11pt]{article}

\usepackage{hyperref}
\hypersetup
{
    colorlinks,
    citecolor=green,
    filecolor=black,
    linkcolor=blue,
    urlcolor=blue
}
\hypersetup{linktocpage}

\usepackage{mathrsfs}

\usepackage{authblk}
\usepackage{blindtext}

\usepackage{graphicx}

\usepackage{amsbsy}
\usepackage{latexsym}
\usepackage{amsfonts}
\usepackage{amssymb}
\usepackage[usenames]{color}
\usepackage{amsmath,amsthm}
\usepackage{enumerate}
\usepackage{authblk}
\usepackage{tikz}
\usetikzlibrary{arrows}
\usepackage{mathdots}
\usepackage{mathtools}
\usepackage{tikz}
\usepackage{cases}
\usepackage{amscd}
\usepackage{amstext}

\usepackage{graphics}
\usepackage{graphicx}
\usepackage{stmaryrd}
\usepackage{dsfont}

\usepackage{color}

\newtheorem{theorem}{Theorem}[section]
\newtheorem{lemma}[theorem]{Lemma}
\newtheorem{prop}[theorem]{Proposition}

\theoremstyle{definition}

\theoremstyle{remark}
\newtheorem{remark}[theorem]{Remark}

\newcommand{\cL}{{\mathcal{L}}}

\newcommand{\cT}{\mathcal{T}}

\newcommand{\dd}{\mathrm{d}}

\newcommand{\cD}{\mathcal{D}}

\newcommand{\cE}{\mathcal{E}}

\newcommand{\Chi}{\raise .3ex
\hbox{\large $\chi$}} 

\newcommand{\T}{{\mathbb{T}}}
\newcommand{\R}{\mathbb{R}}

\newcommand{\N}{\mathbb{N}}
\newcommand{\Z}{\mathbb{Z}}
\newcommand{\C}{\mathbb{C}}

\newcommand{\zz}{\mathbb{Z}}
\newcommand{\re}{\mathbb{R}}
\newcommand{\tor}{\mathbb{T}}

\newcommand\mix{\mathop{\rm mix}}

\usepackage[section]{algorithm}
\usepackage{algorithmicx}
\usepackage{algpseudocode}
\algrenewcommand\algorithmicrequire{\makebox[46pt][l]{\textrm{required:}}}
\algrenewcommand\algorithmicensure{\makebox[46pt][l]{\textrm{output:}}}
\algrenewcommand\algorithmicfunction{\textrm{function}}
\algrenewcommand\algorithmicwhile{\textrm{while}}
\algrenewcommand\algorithmicdo{}
\algrenewcommand\algorithmicend{\textrm{end}}
\algrenewcommand\algorithmicforall{\textrm{for all}}
\algrenewcommand\algorithmicfor{\textrm{for}}
\algrenewcommand\algorithmicrepeat{\textrm{repeat}}
\algrenewcommand\algorithmicuntil{\textrm{until}}
\algrenewcommand\algorithmicif{\textrm{if}}
\algrenewcommand\algorithmicthen{\textrm{then}}
\algrenewcommand\algorithmicelse{\textrm{else}}

\usepackage{todonotes}

\newcommand{\be}{\begin{equation}}
\newcommand{\ee}{\end{equation}}
\newcommand{\beq}{\begin{eqnarray}}
\newcommand{\beqq}{\begin{eqnarray*}}
\newcommand{\eeq}{\end{eqnarray}}
\newcommand{\eeqq}{\end{eqnarray*}}

\numberwithin{equation}{section}

\title{Best $n$-term approximation of diagonal operators and application to function  spaces with mixed smoothness}

\author{Van Kien Nguyen}

\author{Van Dung Nguyen}

\affil{Faculty of Basic Sciences, University of Transport and Communications
	\protect\\	No.3 Cau Giay Street, Lang Thuong Ward, Dong Da District,
	Hanoi, Vietnam
	\protect\\
	Email: kiennv@utc.edu.vn,\
	dungnv@utc.edu.vn
}

\date{\today}
\begin{document}
 
\maketitle
\begin{abstract}  

In this paper we give exact values of the best $n$-term approximation widths of diagonal operators between $\ell_p(\N)$ and $\ell_q(\N)$ with $0<p,q\leq \infty$. The result will be applied to obtain the asymptotic constants of best $n$-term approximation widths of embeddings of function spaces with mixed smoothness by trigonometric system.


\end{abstract}

%&&&&&&&&&&&&&&&&&&&&&&&&&&&&&&&&
%&&&&&&&&&&&&&&&&&&&&&&&&&&&&&&&&

\section{Introduction}
Nowadays, it is well understood that
nonlinear methods of approximation and numerical methods derived from them often produce superior
performance when compared with linear methods. In the last three decades there has been a great success in studying nonlinear approximation which was motivated by numerous  
applications such as numerical analysis, image processing, statistical learning as well as in the design of neural networks. We refer the reader to \cite{DDF.19,DeV98,DeV09} for the development of nonlinear approximation and its application.

In the present paper we concentrate on a particular nonlinear method, the so-called best $n$-term
approximation. 
%The basic idea behind this concept seems quite natural.  We ﬁx a certain system of elements and approximate a given element  by a finite linear combinations of elements contained in this dictionary. 
Our particular interest  is  exact values of best $n$-term approximation of diagonal linear operators from $\ell_p(\N)$ to $\ell_q(\N)$. The exact values of approximation quantities of diagonal operators play an important role in high-dimensional approximation and particularly in studying tractability, see, e.g., \cite{CKS16,KSU14,KSU15,NNS20}. In this paper, the exact values of best $n$-term approximation of diagonal operators will be applied to get the asymptotic constants of best $n$-term approximation of function spaces with mixed smoothness by trigonometric system. 

Let $X$, $Y$ be Banach spaces and $T$  a continuous linear operator from $X$ to $Y$. 
Let $\mathcal{D}$ be a given countable set in $Y$, called dictionary. For given $x\in X$ we  consider the algorithm to approximate $Tx$ by a  finite linear combination of elements contained in this dictionary. 
The
error of this approximation is
\begin{equation*}
\sigma_n(Tx;\mathcal{D}):=\inf_{(a_j)_{j=1}^n\subset \C,\atop (y_j)_{j=1}^n\subset \cD} \Bigg\|Tx-\sum_{j=1}^n a_jy_j\Bigg\|_Y,\qquad n\in \N.
\end{equation*}
We wish to approximate $Tx$ for all $x$ in the closed unit ball of $X$ with respect to the dictionary $\cD$.
This can be measured by the following
benchmark quantity
\begin{equation*}\label{eq-def-sigma}
\sigma_n(T;\mathcal{D}):=\sup_{x\in X,\,\|x\|_X\leq 1}\sigma_n(Tx;\mathcal{D}),\qquad n\in \N.
\end{equation*}
In what follows, we shall call this quantity the best $n$-term approximation width.

Let $\ell_p(\N)$, $0<p\leq \infty$, be the classical complex sequence space with the usual (quasi) norm. 
For $0<p,q\leq \infty$ and positive non-increasing sequence $\lambda=(\lambda_k)_{k\in \N}$,  consider the diagonal linear operator
\begin{equation}\label{eq-diagonal}
T_\lambda:\, (\xi_k)_{k\in \N}  \mapsto (\lambda_k\xi_k)_{k\in \N} 
\end{equation}
 from $\ell_p(\N)$ to $\ell_q(\N)$ and $\cE=\{e_k: k\in \N\}$ where $e_k= (\delta_{k,j})_{j\in \N}$. We are concerned with the exact value of  $\sigma_n(T_\lambda,\cE)$. The first result in this direction was given  by Stepanets \cite{Ste01a}  in the case $p=q$ with the condition $\lim_{k\to \infty}\lambda_k=0$. Later Stepanets generalized his result to the case $0<p\leq q <\infty$, see \cite{Ste01}. Under the same condition $\lim_{k\to \infty}\lambda_k=0$ but by different approach, Gensun and  Lixin \cite{GL06} also obtained exact value of $\sigma_n(T_\lambda,\cE)$ in the case $p=q$.
In this paper we give exact values of  the best $n$-term approximation widths $\sigma_n(T_\lambda,\cE)$, $n\in \N$, in all cases. We also show that the condition $\lim_{k\to \infty}\lambda_k=0$ in the case $0<p<q<\infty$ is not necessary.  Our main result reads as follows. If $0<p<q<\infty$, then we have
\begin{equation*}
 		\sigma_n(T_\lambda,\cE)=
  \dfrac{(n^*-n)^{1/q}}{\big(\sum_{k=1}^{n^*}\lambda_k^{-p}\big)^{1/p}},   
\end{equation*}
where $n^*$ is the smallest integer $m>n$ such that 
$$
\dfrac{(m-n)^{1/q}}{\big(\sum_{k=1}^m\lambda_k^{-p}\big)^{1/p}} \geq  \dfrac{(m+1-n)^{1/q}}{\big(\sum_{k=1}^{m+1}\lambda_k^{-p}\big)^{1/p}}.
$$  
In the case $0<q<p<\infty$ and the series $\sum_{k=1}^\infty \lambda_k^{pq/(p-q)}$ converges, we get
\begin{equation*}
	\sigma_n(T_\lambda,\cE)=
	\Bigg(\dfrac{(n_*-n)^{\frac{p}{p-q}}}{\big(\sum_{k=1}^{n_*}\lambda_k^{-p}\big)^{\frac{q}{p-q}}} + \sum_{k=n_*+1}^\infty \lambda_k^{\frac{pq}{p-q}}\Bigg)^{\frac{1}{q}-\frac{1}{p}},
\end{equation*}
where $n_*$ is the largest integer $m>n$ such that 
\begin{equation*}
(m-n)\lambda_m^{-p} \leq \sum_{k=1}^m \lambda_k^{-p}.
\end{equation*}
The limiting cases $p=q$ or $p=\infty$ and/or $q=\infty$ are also obtained, see Theorem \ref{prop-sequence}.

The above results will be applied to study best $n$-term approximation of embedding of function spaces with mixed smoothness by trigonometric system  $\cT^d:=\{e^{{\rm i}k x}:\, k \in \Z^d\}$ on the torus $\T^d$ of  dimension $d$.  
 Our motivation stems from high-dimensional 
approximation which has been the object of an intensive study recently. 
In many high-dimensional 
approximation problems when the high-dimensional signals or functions have appropriate mixed smoothness, one can apply efficiently approximation methods and sampling algorithms constructed on sparse grids to obtain tractability for algorithms or numerical methods. We refer the reader to the monographs \cite{NoWo08,NoWo10}  for concepts of computation complexity and results on high dimensional problems. Survey on various aspects of  high-dimensional approximation of functions having  mixed smoothness can be found in the recent book  \cite{DTU18B}.

There has been a numerous papers working on  best $n$-term approximation of embeddings of function spaces with mixed smoothness by different dictionaries. For instance, Bazarkhanov \cite{Baz16}, Dinh D\~ung \cite{Dung98a,Dung00,Dung01}, Kashin and Temlyakov \cite{KT94},  Romanyuk 
\cite{Roma03,Roma06}, Romanyuk and Romanyuk \cite{RoRo10}, Temlyakov \cite{Tem86b,Tem00,Tem86a,Tem15} worked on trigonometric system; Hansen and Sickel \cite{BS17,Dung01b,HaSi10,HaSi12} on wavelet system. For some recent contributions in this direction we refer to 
\cite{BeTo19,Byr18,Tem17,TU21}. Historical
comments and further references for studies of best $n$-term approximation of function spaces with mixed smoothness can be found in the two recent books \cite[Chapter 7]{DTU18B} and \cite[Chapter 9]{Tem18B}. Let us emphasize  here that all the above mentioned papers worked only on the asymptotic order of best $n$-term approximation widths of embeddings of function spaces with mixed smoothness. The asymptotic constants and pre-asymptotic estimates of this quantity have not been considered.

Let  $0<s<\infty$ and $0<r\leq \infty$. This paper considers the best $n$-term approximation of embedding of Sobolev space with mixed smoothness $H^{s,r}_{\mix}(\T^d)$ on the torus $\T^d$ into either $L_2(\T^d)$ or Wiener space $\mathcal{A}(\T^d)$. 
In this context we will not only investigate the optimal order of the decay of the best $n$-term approximation widths   but we will determine 
the asymptotic constant as well. This sheds some light not only on the dependence on $n$, but also on the dependence on $s,r$ 
and in particular on $d$.  We have 
 	\begin{equation*} 
 		\lim\limits_{n\to \infty} \frac{\sigma_n\big(id: H^{s,r}_{\mix}(\T^d)\to L_2(\T^d),\cT^d\big)}{n^{-s}(\ln n)^{s(d-1)}}=\frac{s^s}{(s+\frac{1}{2})^s} \bigg( \frac{2^d}{(d-1)!}\bigg)^s
 	\end{equation*}
and if $s>1/2$
\begin{equation*} 
	\lim\limits_{n\to \infty} \frac{\sigma_n\big(id: H^{s,r}_{\mix}(\T^d)\to \mathcal{A}(\T^d),\cT^d\big)}{n^{-s+\frac{1}{2}}(\ln n)^{s(d-1)}}=\bigg(\frac{s}{s+\frac{1}{2}}\bigg)^s\bigg(\frac{1}{ s-\frac{1}{2}}\bigg)^{\frac{1}{2}}\bigg( \frac{2^d}{(d-1)!}\bigg)^s\,.
\end{equation*}

In this paper we also obtain the asymptotic constants of best $n$-term approximation widths of embeddings of Sobolev spaces with mixed smoothness $H^{s,2}_{\mix}(\T^d)$ into  the  energy norm space $H^1(\T^d)$. Those embeddings are of particular importance with respect to the numerical solution of the Poisson equation, see \cite{BuGr04}. In this case, with $s>1$ we get
	\begin{equation*}
		\lim\limits_{n\to \infty} \frac{ \sigma_n\big(id: H_{\mix}^{s,2}(\T^d) \to H^1(\T^d),\cT^d\big)}{n^{-s+1}}= \frac{(s-1)^{s-1}}{(s-\frac{1}{2})^{s-1}} (2d)^{s-1}(2S+1)^{(s-1)(d-1)}, 
	\end{equation*}
where 
$$
S:=\sum_{k=1}^{+\infty}\frac{1}{(k^2+1)^{\frac{s}{2(s-1)}}}\, .
$$

%&&&&&&&&&&&&&&&&&&&&&&&&&&&&&&&&&&&&&&&&&&&&&&&&&&&&&&&&&&&&&&&&&&&&&&&&&&&&&&&&&&&&&&&&&&&&&&&&&&&&&&&&&&&&&&&
%&&&&&&&&&&&&&&&&&&&&&&&&&&&&&&&&&&&&&&&&&&&&&&&&&&&&&&&&&&&&&&&&&&&&&&&&&&&&&&&&&&&&&&&&&&&&&&&&&&&&&&&&&&&&&&&

The paper is organized as follows. 
In Section \ref{sec-diagonal} we collect some properties of best $n$-term approximation widths and give  exact values of best $n$-term approximation widths of diagonal operators. 
The next Section \ref{sec-F-Class} is devoted to the study of 
asymptotic constants of best $n$-term approximation widths of embeddings of weighted  classes  $F_{\omega,p}(\T^d)$.
These results will be used in final Section \ref{sec-mixed}, where we deal with the particular family  of weights  
associated to  function spaces of dominating mixed smoothness.

%&&&&&&&&&&&&&&&&&&&&&&&&&&&&&&&&
%&&&&&&&&&&&&&&&&&&&&&&&&&&&&&&&&

\subsection*{Notation}

%&&&&&&&&&&&&&&&&&&&&&&&&&&&&&&&&&&
%&&&&&&&&&&&&&&&&&&&&&&&&&&&&&&&&&&

As usual, $\N$ denotes the natural numbers, $\N_0$ the non-negative integers,
$\zz$ the integers,
$\re$ the real numbers, and $\C$ the complex numbers. We denote by $\T$  
 the torus, represented by the interval $[0,2\pi]$, where
the end points of the interval are identified.
For a real number $a$ we denote by $\lfloor a \rfloor$ the greatest integer not larger than $a$.
The letter $d$ is always reserved for the dimension in $\N^d$, $\Z^d$, $\re^d$, $\C^d$, and $\tor^d$.
%For $0<p\leq \infty$ and $x=(x_1,\ldots,x_d)\in \re^d$ we denote $|x|_p = \big(\sum_{i=1}^d |x_i|^p\big)^{1/p}$ with the usual modification for $p=\infty$. 
%If $\alpha=(\alpha_1,\ldots,\alpha_d)\in \N_0^d$ and $x \in \C^d$ we use
%$x^\alpha:=\prod_{i=1}^d x_i^{\alpha_i}$ with the convention $0^0:=1$. 
For two Banach spaces $X$ and $Y$, $\cL(X,Y)$ denotes the set of  continuous linear operators from $X$ to $Y$. If $(a_n)_{n\in \N}$ and $(b_n)_{n\in \N}$ are two sequences, the symbol $a_n\sim b_n, n\to \infty$, indicates that $\lim_{n\to \infty}\frac{a_n}{b_n}=1$.
The equivalence $a_n\asymp b_n$ means that there are constants $0<c_1\le c_2<\infty$ such that
$c_1\, a_n\le b_n\le c_2\, a_n$ for all $n\in\N$.

%&&&&&&&&&&&&&&&&&&&&&&&&&&&&&&&&&&&&&&&&&
%&&&&&&&&&&&&&&&&&&&&&&&&&&&&&&&&&&&&&&&&&

%&&&&&&&&&&&&&&&&&&&&&&&&&&&&&&&&&&&&&&&&&
%&&&&&&&&&&&&&&&&&&&&&&&&&&&&&&&&&&&&&&&&&

\section{Best $n$-term approximation widths of diagonal operators}\label{sec-diagonal}

%&&&&&&&&&&&&&&&&&&&&&&&&&&&&&&&&&&&&&&&&&
This section is devoted to give exact values of the best $n$-term approximation widths of the diagonal operator defined in \eqref{eq-diagonal}. Let $X,Y$ be Banach spaces, $T\in \mathcal{L}(X,Y)$, and $\mathcal{D}\subset Y$  a dictionary. By definition, it is clear that $(\sigma_n(T,\cD))_{n\in \N}$ is a non-increasing sequence. If  $W, Z$ are Banach spaces and $A\in \cL(W,X)$, $B\in \cL(Y,Z)$ then we have
\begin{equation}\label{eq-ideal}
			\sigma_n(BTA,B(\cD)) \leq \|B\| \cdot \sigma_n(T,\cD) \cdot \|A\|.
		\end{equation}
A proof of this fact can be found in 
\cite[Lemma 6.1]{Byr18}. 
For further properties of the  best $n$-term approximation widths such as additivity, interpolation we refer the reader to  \cite{Han10,Byr18,TU21}. In fact the best $n$-term approximation widths belong to the notion of pseudo $s$-numbers
introduced by Pietsch, see \cite{TU21}. 

%&&&&&&&&&&&&&&&&&&&&&&&&&&&&&&&&&&&&&&&&&

Let $
T_\lambda$
be the diagonal operator from $\ell_p(\N)$ to $\ell_q(\N)$ defined in \eqref{eq-diagonal}. By definition we have
\begin{equation}\label{eq-def01}
\sigma_n(T_\lambda,\cE)= 
\begin{cases}
	\sup\limits_{(\xi_k)_{k\in \N}\in B_p}\inf\limits_{\Gamma_n} \bigg(\sum\limits_{k\not\in \Gamma_n}|\lambda_k \xi_k|^q\bigg)^{1/q} &\text{if } 0< q<\infty
	\\
	\sup\limits_{(\xi_k)_{k\in \N}\in B_p}\inf\limits_{\Gamma_n}  \sup\limits_{k\not\in \Gamma_n}|\lambda_k \xi_k| &\text{if }   q=\infty,
\end{cases}
\end{equation}
where $B_p$ is the closed unit ball of $\ell_p(\N)$ and $\Gamma_n$ is an arbitrary subset of $\N$ with $n$ elements. In the following we give exact value of $\sigma_n(T_\lambda,\cE)$ with $0<p,q\leq \infty$. The proof is   mainly based on the exact values of  $n$-term approximation widths of the diagonal operators from $\ell_p^M$ to $\ell_q^M$ obtained  by Gao in \cite{Gao10}. Here $\ell_p^M$ stands for $\C^M$ equipped with the usual norm $\|\cdot\|_{\ell_p^M}$. 
For a vector $\lambda=(\lambda_j)_{j=1}^M$ with $\lambda_1\geq \lambda_2\geq \ldots \geq \lambda_M > 0$ the diagonal operator $T^M_\lambda$ from $\ell_p^M$ to $\ell_q^M$ is defined by $(\xi_j)_{j=1}^M\mapsto (\lambda_j\xi_j)_{j=1}^M$. Let  $\cE_M=\{e_1,\ldots,e_M\}$ be  the standard basis of $\R^M$. If $n\in \N$ and $n\leq M$ we have 
\begin{equation}\label{eq-def-sigmalM}
\sigma_n(T_\lambda^M,\cE_M)=
	\sup\limits_{(\xi_k)_{k=1}^M\in B_p^M}\inf\limits_{\Gamma_n^M} \bigg(\sum\limits_{k\not\in \Gamma_n^M}|\lambda_k \xi_k|^q\bigg)^{1/q}, \quad 0< q<\infty,
\end{equation}
where $B_p^M$ is the closed unit ball of $\ell_p^M$ and $\Gamma_n^M$ is an arbitrary subset of $\{1,\ldots,M\}$ with $n$ elements. For $\lambda=(\lambda_j)_{j\in \N}$ we define $T^M_\lambda:=T^M_{\tilde{\lambda}}$ where $\tilde{\lambda}=(\lambda_j)_{j=1}^M$.
When $q=\infty$ the summation in \eqref{eq-def-sigmalM} is replaced by supremum.

 Note that if $\lambda=(\lambda_j)_{j\in \N}$ satisfying $\lambda_1\geq \lambda_2\geq \ldots \geq \lambda_M >0$ and $\lambda_j=0$ for $j\geq M+1$, then
$$
 \sigma_n(T_\lambda,\cE) = \sigma_n(T^M_\lambda,\cE_M),\quad n\in \N
$$
which were obtained in \cite{Gao10}. Therefore, we only consider the operator $T_\lambda$ where  $\lambda=(\lambda_j)_{j\in \N}$ is a positive sequence.
Our main result in this section reads as follows.
\begin{theorem}\label{prop-sequence}
Let $0<p,q\leq \infty$ and $\lambda=(\lambda_k)_{k\in \N}$ be a positive non-increasing sequence. Let  $T_\lambda$ be defined in \eqref{eq-diagonal} and $n\in \N$. 
\begin{description}
    \item[(i)] If $0<p\leq q<\infty$ we have
\begin{equation*}
 		\sigma_n(T_\lambda,\cE)=
 	\sup_{m> n} \dfrac{(m-n)^{1/q}}{\big(\sum_{k=1}^m\lambda_k^{-p}\big)^{1/p}}.   
\end{equation*}
Moreover, if either  $0<p<q<\infty$ or $\lim\limits_{k\to \infty}\lambda_k=0$ then
\begin{equation*}
 		\sigma_n(T_\lambda,\cE)=
  \dfrac{(n^*-n)^{1/q}}{\big(\sum_{k=1}^{n^*}\lambda_k^{-p}\big)^{1/p}}\,,   
\end{equation*}
where $n^*$ is the smallest integer $m>n$ such that 
$$
\dfrac{(m-n)^{1/q}}{\big(\sum_{k=1}^m\lambda_k^{-p}\big)^{1/p}} \geq \dfrac{(m+1-n)^{1/q}}{\big(\sum_{k=1}^{m+1}\lambda_k^{-p}\big)^{1/p}}.
$$
\item[(ii)] If $0<q<p<\infty$ and the series $\sum_{k=1}^\infty \lambda_k^{pq/(p-q)}$ converges, then we have
\begin{equation}\label{eq-q<p}
	\sigma_n(T_\lambda,\cE)=
	\Bigg(\dfrac{(n_*-n)^{\frac{p}{p-q}}}{\big(\sum_{k=1}^{n_*}\lambda_k^{-p}\big)^{\frac{q}{p-q}}} + \sum_{k=n_*+1}^\infty \lambda_k^{\frac{pq}{p-q}}\Bigg)^{\frac{1}{q}-\frac{1}{p}}\,,
\end{equation}
where $n_*$ is the largest integer $m>n$ such that 
\begin{equation}\label{eq-n-*}
(m-n)\lambda_m^{-p} \leq \sum_{k=1}^m \lambda_k^{-p}.
\end{equation}
\item[(iii)] If $0<p<q=\infty$ then
\begin{equation*}
 		\sigma_n(T_\lambda,\cE)=
 	   \bigg(\sum_{k=1}^{n+1}\lambda_k^{-p}\bigg)^{-1/p}.   
\end{equation*}
\item[(iv)] If $0<q<p=\infty$ and the series $\sum_{k=1}^\infty \lambda_k^{q}$ converges then
\begin{equation*}
	\sigma_n(T_\lambda,\cE)=
	\Bigg( \sum_{k=n+1}^\infty \lambda_k^{q}\Bigg)^{1/q}.
\end{equation*}
\item[(v)] If $p=q=\infty$  then
\begin{equation*}
	\sigma_n(T_\lambda,\cE)=\lambda_{n+1}.
\end{equation*}
\end{description} 
\end{theorem}
As mentioned in Introduction, the exact values of $\sigma_n(T_\lambda,\cE)$, $n\in \N$, in the case $0<p\leq q\leq \infty$ were obtained in \cite{Ste01a,Ste01,GL06} under the condition $\lim_{k\to \infty}\lambda_k=0$. To prove the above theorem, we need some auxiliary results.
\begin{lemma} \label{lem-n*}
Let $0<p<q<\infty$ and $(\lambda_k)_{k=1}^\infty$ be a positive non-increasing sequence. Then for $n\in \N$, there is $n^*=n^*(n)\in \N$ such that
\begin{equation}\label{eq-n*-k2}
     	\sup_{m> n} \dfrac{(m-n)^{1/q}}{\big(\sum_{k=1}^m\lambda_k^{-p}\big)^{1/p}} = \dfrac{(n^*-n)^{1/q}}{\big(\sum_{k=1}^{n^*}\lambda_k^{-p}\big)^{1/p}}.
\end{equation}
Moreover, $n^*$ can be chosen as the smallest integer $m>n$ such that 
$$
\dfrac{(m-n)^{1/q}}{\big(\sum_{k=1}^m\lambda_k^{-p}\big)^{1/p}} \geq \dfrac{(m+1-n)^{1/q}}{\big(\sum_{k=1}^{m+1}\lambda_k^{-p}\big)^{1/p}}.
$$
\end{lemma}
\begin{proof}
We first show that $n^*$ exists. The case $\lim_{k\to \infty}\lambda_k=0$ was already considered in \cite{Ste01}. We prove  the case $\lim_{k\to \infty}\lambda_k=K>0$.  We will show that there exists $n_0\in \N$ such that 
\begin{equation}\label{eq-n*-k1}
 \dfrac{(m-n)^{1/q}}{\big(\sum_{k=1}^m\lambda_k^{-p}\big)^{1/p}} < \dfrac{(2n-n)^{1/q}}{\big(\sum_{k=1}^{2n}\lambda_k^{-p}\big)^{1/p}}
\end{equation}
for $m>n_0$ and as a consequence 
we obtain \eqref{eq-n*-k2}
for some $n^*\in \{n+1,\ldots,n_0\}$. 
Observe that if $m\in \{jn+1,\ldots, (j+1)n \}$ for some $j\in \N$ we have
\begin{equation*}
    \begin{aligned}
 \dfrac{(m-n)^{1/q}}{\big(\sum_{k=1}^m\lambda_k^{-p}\big)^{1/p}} < \dfrac{(jn)^{1/q}}{\big(\sum_{k=1}^{jn}\lambda_k^{-p}\big)^{1/p}} \leq \dfrac{(jn)^{1/q}}{\lambda_1^{-1} (jn)^{1/p}}.
    \end{aligned}
\end{equation*}
We also have
\begin{equation*}
    \begin{aligned}
   \dfrac{n^{1/q}}{K^{-1}  (2n)  ^{1/p}} \le \dfrac{(2n-n)^{1/q}}{\big(\sum_{k=1}^{2n}\lambda_k^{-p}\big)^{1/p}}.
    \end{aligned}
\end{equation*}
Therefore 
$$
\sup_{m\in \{jn+1,\ldots, (j+1)n \}}\dfrac{(m-n)^{1/q}}{\big(\sum_{k=1}^m\lambda_k^{-p}\big)^{1/p}} < \dfrac{(2n-n)^{1/q}}{\big(\sum_{k=1}^{2n}\lambda_k^{-p}\big)^{1/p}}
$$
if
\begin{equation*}
    \begin{aligned}
 \dfrac{(jn)^{1/q}}{\lambda_1^{-1} (jn)^{1/p}} <   \dfrac{n^{1/q}}{K^{-1}  (2n)  ^{1/p}} 
\quad \Longleftrightarrow \quad
j  >  \bigg(\dfrac{\lambda_1  2^{1/p}}{K} \bigg)^{\frac{pq}{q-p}} .
    \end{aligned}
\end{equation*}
Denoting
$
n_0 =\big\lceil \big(\frac{\lambda_1  2^{1/p}}{K} \big)^{ \frac{pq}{q-p}} \big\rceil n
$
we obtain \eqref{eq-n*-k1} for $m>n_0$ and \eqref{eq-n*-k2} follows.

We turn to the second statement.
Assume 
\begin{equation}\label{eq-assumption-1}
\dfrac{m_0-n}{\big(\sum_{k=1}^{m_0}\lambda_k^{-p}\big)^{q/p}}\geq \dfrac{m_0+1-n}{\big(\sum_{k=1}^{m_0+1}\lambda_k^{-p}\big)^{q/p}}
\end{equation}
for some $m_0>n$, $m_0\in \N$. Such $m_0$ exists by the first statement. We will prove
\begin{equation}\label{eq-objective0}
      \dfrac{m_0+1-n}{\big(\sum_{k=1}^{m_0+1}\lambda_k^{-p}\big)^{q/p}}
  \geq   \dfrac{m_0+2-n}{\big(\sum_{k=1}^{m_0+2}\lambda_k^{-p}\big)^{q/p}}
\end{equation}
by showing that
\begin{equation}\label{eq-objective}
      \dfrac{m_0+1-n}{\big(\sum_{k=1}^{m_0+1}\lambda_k^{-p}\big)^{q/p}}
  \geq   \dfrac{m_0+2-n}{\big(\sum_{k=1}^{m_0+1}\lambda_k^{-p}+\lambda_{m_0+1}^{-p}\big)^{q/p}}\,.
\end{equation}
Putting $A=\sum_{k=1}^{m_0}\lambda_k^{-p}$ and $a=\lambda_{m_0+1}^{-p}$ we consider the function
$
g(h)=  \frac{m_0-n+h}{(A+ha)^{q/p}}\,.
$ We have
$$
g'(h) = \dfrac{(A+ha)-\frac{q}{p}a(m_0-n+h)}{(A+ha)^{\frac{q}{p}+1}} 
$$
and
$g'(h) \le  0$ if 
$$
(A+ha)-\frac{q}{p}a(m_0-n+h) \le  0 \ \ \Longleftrightarrow \ \ h\geq \dfrac{A-\frac{q}{p}a(m_0-n)}{a(\frac{q}{p}-1)}. 
$$
Assume   
\begin{equation}\label{eq-assumption1}
\dfrac{A-\frac{q}{p}a(m_0-n)}{a(\frac{q}{p}-1)} \ge 1 \ \ \Longleftrightarrow \ \ \frac{a}{A} \le \frac{1}{\frac{q}{p}(m_0-n) +\frac{q}{p} -1}\,.
\end{equation}
Observe that the condition \eqref{eq-assumption-1} implies 
$\big(
1+\frac{a}{A}
\big)^{q/p} \geq 1+\frac{1}{m_0-n}.
$
From this and \eqref{eq-assumption1} we get
$$
\bigg(
1+\frac{1}{\frac{q}{p}(m_0-n) +\frac{q}{p} -1}
\bigg)^{q/p} \geq 1+\frac{1}{m_0-n}.
$$
But this is a contradiction since  $\varphi(t)=(1+\frac{1}{t(m_0-n)+t-1})^t$ is a strictly decreasing function on $[1,+\infty)$ and $\varphi(1) = 1+\frac{1}{m_0-n}$. Consequently, $g$ is  decreasing on $[1,+\infty)$. This proves \eqref{eq-objective} and \eqref{eq-objective0} follows.
\end{proof}
\begin{lemma}\label{lem-elementary}
Let $(\delta_k)_{k\in \N}$ be a positive increasing  sequence and $\lim\limits_{k\to \infty}\delta_k=+\infty$. Let $n\in \N$ and $n_*=n_*(n)$ be the largest integer $m>n$ such that
\begin{equation}\label{eq-n8}
	(m-n)\delta_m \leq \sum_{k=1}^m \delta_k.
\end{equation}
Then $n_*$ is finite and for any $m\in \{n+1,\ldots,n_*\}$ the inequality \eqref{eq-n8} holds true.
\end{lemma}
\begin{proof} First of all, observe that $m=n+1$ satisfies \eqref{eq-n8}. If $m\geq n+1$ and $m$ satisfies \eqref{eq-n8}  we can write
	\begin{align*}
		m & \leq n+ \frac{\delta_1+\ldots + \delta_{n+1}}{\delta_m} +   \frac{ \delta_{n+2}+\ldots+\delta_{m}}{\delta_m}
		\\
		& \leq n+ \frac{\delta_1+\ldots + \delta_{n+1}}{\delta_m} + m-n-1 
=m-1 + \frac{\delta_1+\ldots + \delta_{n+1}}{\delta_m}.
	\end{align*}
	Observe that the term $\frac{\delta_1+\ldots + \delta_{n+1}}{\delta_m}$ tends to zero when $m\to \infty$. Consequently $n_*$ is finite.
 Assume  
\begin{equation}\label{eq-n0-1}
	(m_0-n)\delta_{m_0} > \sum_{k=1}^{m_0} \delta_k
\end{equation}
for some $m_0\in \N$, $m_0>n$.
Then we have
\begin{equation*}
	\begin{aligned}
		m_0+1 > \frac{1}{\delta_{m_0}}\sum_{k=1}^{m_0} \delta_k +n+1 > \frac{1}{\delta_{m_0+1}}\sum_{k=1}^{m_0} \delta_k+1+n=\frac{1}{\delta_{m_0+1}}\sum_{k=1}^{m_0+1} \delta_k+n.
	\end{aligned}
\end{equation*}
This shows that the inequality \eqref{eq-n0-1} is satisfied with $m_0$ being replaced by $m_0+1$ and therefore is satisfied with any  $m>m_0$, $m\in \N$. As a consequence we conclude that for any $m\in \{n+1,\ldots,n_*\}$ the inequality \eqref{eq-n8} holds true.
\end{proof}
We are now in position to prove Theorem \ref{prop-sequence}.
\begin{proof}
{\it Step 1.} Proof of (i). Given $\varepsilon>0$. For any $(\xi_k)_{k\in \N}\in B_p$ we take $M\in \N$ (depending on $(\xi_k)_{k\in \N}$) such that
$$
 \sum_{k=M+1}^\infty |\xi_k|^p<\varepsilon.
$$
Then we have
\begin{equation}\label{eq-sigma-M}
    \begin{aligned}
\inf_{\Gamma_n}\bigg(\sum_{k\not \in \Gamma_n }|\lambda_k \xi_k|^q\bigg) 
& \leq \inf_{\Gamma_n^M}\bigg(\sum_{k\not \in \Gamma_n^M }|\lambda_k \xi_k|^q\bigg) 
 = \inf_{\Gamma_n^M}\bigg(\sum_{k\in\{1,\ldots,M\} \backslash \Gamma_n^M }|\lambda_k \xi_k|^q\bigg)+\sum_{k=
	M+1}^\infty|\lambda_k \xi_k|^q.
\end{aligned}
\end{equation}
The first term on the right side can be estimated as follows
\begin{equation}\label{eq-kien01}
\inf_{\Gamma_n^M}\bigg(\sum_{k\in\{1,\ldots,M\} \backslash \Gamma_n^M }|\lambda_k \xi_k|^q\bigg)\leq   \sup_{(\gamma_k)_{k=1}^M\in B_p^M} \inf_{\Gamma_n^M}\bigg(\sum_{k\in\{1,\ldots,M\} \backslash \Gamma_n^M }|\lambda_k \gamma_k|^q\bigg) = \sigma_n(T_\lambda^M,\cE_M)^q,
\end{equation}
see \eqref{eq-def-sigmalM}. It has been proved in \cite{Gao10} that
\begin{equation}\label{eq-Gao1}
\sigma_n(T^M_\lambda,\cE_M) = \sup_{n<m\leq M} \dfrac{(m-n)^{1/q}}{\big(\sum_{k=1}^m\lambda_k^{-p}\big)^{1/p}}.
\end{equation}
Hence we get
\begin{equation*}
    \begin{aligned}
\inf_{\Gamma_n^M}\bigg(\sum_{k\in\{1,\ldots,M\} \backslash \Gamma_n^M }|\lambda_k \xi_k|^q\bigg)
&  \leq	\sup_{n<m\leq M} \dfrac{m-n}{\big(\sum_{k=1}^m\lambda_k^{-p}\big)^{q/p}}
 \leq \sup_{n<m} \dfrac{m-n}{\big(\sum_{k=1}^m\lambda_k^{-p}\big)^{q/p}}. 
\end{aligned}
\end{equation*}
Since $0<p\leq q\leq \infty$, for the second term we have
\begin{equation*}
    \begin{aligned}
\sum_{k=M+1}^{\infty}|\lambda_k\xi_k|^q \leq \bigg(\sum_{k=
	M+1}^\infty|\lambda_k	\xi_k|^p\bigg)^{q/p}
& \leq \sup_{k>M}|\lambda_k|^q\bigg(\sum_{k=
	M+1}^\infty|	\xi_k|^p\bigg)^{q/p} \leq  \lambda_1^q\varepsilon^{q/p}   .
    \end{aligned}
\end{equation*}
Consequently we obtain
\begin{equation*}
    \begin{aligned}
\inf_{\Gamma_n}\bigg(\sum_{k\not \in \Gamma_n }|\lambda_k \xi_k|^q\bigg)
& \leq 	\sup_{m> n} \dfrac{m-n}{\big(\sum_{k=1}^m\lambda_k^{-p}\big)^{q/p}}+ \lambda_1^q \varepsilon^{q/p}  .
\end{aligned}
\end{equation*}
Observe that the right-hand side is independent of $(\xi_k)_{k\in \N}\in B_p$ and $\varepsilon>0$ is arbitrarily small. In view of \eqref{eq-def01} we obtain the upper bound. 

We now give a proof for the lower bound. Take $M\in \N$ arbitrarily large and consider the following diagram 
		\begin{equation*}
		\begin{CD}
			\ell_p^M  @ > T_\lambda^M >> \ell_q^M \\
			@VV J V @AA Q A\\
			\ell_p(\N) @ > T_\lambda >> \ell_q(\N) \, ,
		\end{CD}
	\end{equation*}
where
\begin{equation*}\label{eq-TM}
\begin{aligned}
&J(\xi_1,\ldots,\xi_M) =(\xi_1,\ldots,\xi_M,0,0,\ldots)
\\
&Q(\xi_1,\ldots,\xi_M,\xi_{M+1},\ldots) = (\xi_1,\ldots,\xi_M).
\end{aligned}	
\end{equation*}
We have $T_\lambda^M=QT_\lambda J$ and $\|J\|=\|Q\|=1$ which by property \eqref{eq-ideal} implies 
$$
\sigma_n(T^M_\lambda,\cE_M)  \leq \|J\| \cdot \sigma_n(T_\lambda,\cE)\cdot \|Q\| = \sigma_n(T_\lambda,\cE).
$$
Using \eqref{eq-Gao1} again we deduce 
$$
 \sup_{n<m\leq M} \dfrac{(m-n)^{1/q}}{\big(\sum_{k=1}^m\lambda_k^{-p}\big)^{1/p}} \leq \sigma_n(T_\lambda,\cE).
$$
Since $M$ is arbitrarily large, we obtain the lower bound. The second statement follows from Lemma \ref{lem-n*}.
\\
{\it Step 2.} Proof of (ii). First note that $n_*$ is finite by Lemma \ref{lem-elementary}. Let $\varepsilon >0$. We choose $M>n_*$ such that
$$
\bigg(\sum_{k=M+1}^\infty \lambda_k^{\frac{pq}{p-q}}\bigg)^{\frac{p-q}{p}} <\varepsilon.
$$ 
For $(\xi_k)_{k\in \N}\in B_p$, we use the estimate \eqref{eq-sigma-M}.
Applying H\"older's inequality we get
\begin{align*}
\sum_{k=
	M+1}^\infty|\lambda_k \xi_k|^q \leq \bigg(\sum_{k=M+1}^\infty \lambda_k^{\frac{pq}{p-q}}\bigg)^{\frac{p-q}{p}} \bigg( \sum_{k=M+1}^\infty |\xi_k|^{p}\bigg)^{{\frac{q}{p}}} <\varepsilon
\end{align*}
which by \eqref{eq-sigma-M} implies
\begin{equation*}
	\begin{aligned}
 \inf_{\Gamma_n}\bigg(\sum_{k\not \in \Gamma_n }|\lambda_k \xi_k|^q\bigg) 
		& \leq  \inf_{\Gamma_n^M}\bigg(\sum_{k\in\{1,\ldots,M\} \backslash \Gamma_n^M }|\lambda_k \xi_k|^q\bigg) +\varepsilon
 = \sigma_n(T^M_\lambda,\cE_M)^q + \varepsilon,
	\end{aligned}
\end{equation*}
see \eqref{eq-kien01}. Using the result in \cite{Gao10} for the case $0<q<p<\infty$
\begin{equation*}
	\begin{aligned}
\sigma_n(T^M_\lambda,\cE_M)
& =
\Bigg(\dfrac{(n_*-n)^{\frac{p}{p-q}}}{\big(\sum_{k=1}^{n_*}\lambda_k^{-p}\big)^{\frac{q}{p-q}}} + \sum_{k=n_*+1}^M \lambda_k^{\frac{pq}{p-q}}\Bigg)^{\frac{p-q}{pq}}  
\leq
\Bigg(\dfrac{(n_*-n)^{\frac{p}{p-q}}}{\big(\sum_{k=1}^{n_*}\lambda_k^{-p}\big)^{\frac{q}{p-q}}} + \sum_{k=n_*+1}^\infty \lambda_k^{\frac{pq}{p-q}}\Bigg)^{\frac{p-q}{pq}}  
	\end{aligned}
\end{equation*} 
and following the argument as in Step 1 we obtain the upper bound. The lower bound is carried out similarly as Step 1 with $M>n_*$. The other cases are proved similarly with a slight modification.  
\end{proof}

%&&&&&&&&&&&&&&&&&&&&&&&&&&&&&&&&&&&&&&&&&
%&&&&&&&&&&&&&&&&&&&&&&&&&&&&&&&&&&&&&&&&&
\section{Best $n$-term approximation of function classes  $F_{\omega,p}(\T^d)$}\label{sec-F-Class}

%&&&&&&&&&&&&&&&&&&&&&&&&&&&&&&&&&&&&&&&&&
%&&&&&&&&&&&&&&&&&&&&&&&&&&&&&&&&&&&&&&&&&
Let $\T^d$ be the $d-$dimensional torus. We equip $\T^d$ with the probability measure $(2\pi)^{-d}\dd x$. 
In this section we study the asymptotic constants of best $n$-term approximation widths of embeddings of the weighted function classes $F_{\omega,p}(\T^d)$ by trigonometric system $\cT^d$. 
For a function $f\in L_1(\T^d)$, its Fourier coefficients are defined as
\begin{equation*}
	\hat{f}(k):=(2\pi)^{-d}\int_{\T^d} f(x) e^{-{\rm i} k x} \dd x,\qquad k\in\Z^d.
\end{equation*}
Hence, it holds for any $f\in L_2(\T^d)$ that
\begin{equation*}
	\|f\|_{L_2(\T^d)}^2= (2\pi)^{-d}\int_{\T^d}|f(x)|^2  \dd x = \sum_{k \in \Z^d} |\hat{f}(k)|^2 \,.
\end{equation*}

Let $\omega=(\omega(k))_{k\in \Z^d}$  be a sequence of positive   numbers. Those sequences we will call a weight in what follows.
For $0< p \leq \infty$ we introduce the  class 
$F_{\omega,p}(\T^d)$ as the collection of all  functions $f\in L_1(\T^d)$ such that 
\[
\|f\|_{F_{\omega,p}(\T^d)}:= \bigg( \sum_{k\in \Z^d} |\omega(k) \hat{f}(k)|^p\bigg)^{1/p}< \infty \, .
\]
When $\omega(k)= 1$ for all $k\in \Z^d$ we use the notation $F_p(\T^d)$ instead of $F_{\omega,p}(\T^d)$. In this case we get back the space $L_2(\T^d)$ when $p=2$ and the classical  Wiener algebra $\mathcal{A}(\T^d)$ when $p=1$. 

We suppose that  
\begin{equation}\label{eq-condition}
\lim_{|k_1|+\ldots+|k_d|\to \infty} \omega(k) =+\infty,\quad k=(k_1,\ldots,k_d).
\end{equation}
In what follows we denote the  non-increasing rearrangement of the sequence $(1/\omega(k))_{k\in \Z^d}$ 
  by  $\lambda=(\lambda_n)_{n\in \N}$. 
Observe that   $id: F_{\omega,2}(\T^d) \to  L_2(\T^d)$ is compact if and only if
 $\lim_{n\to \infty}\, \lambda_n =0.$
In fact we have
$$\lambda_n = a_n(id:F_{\omega,2}(\T^d) \to  L_2(\T^d) ),$$
where $a_n(id:F_{\omega,2}(\T^d) \to  L_2(\T^d) )$ is the  $n$-th approximation number (linear width) of the operator $id: F_{\omega,2}(\T^d) \to  L_2(\T^d)$, see  \cite{KSU15}. Recall that for two Banach spaces $X$, $Y$ and $T\in \cL(X,Y)$, the $n$-th approximation number of $T$ is defined as
$$ a_n(T):=\inf\big\{\|T-A:X\to Y\|: \ A\in \mathcal L(X,Y),\ \ \text{rank} (A)<n\big\}\, , \quad n \in \N\, . $$
Basic properties of this quantity can be found in \cite{Pie80B,Pie87B}.

We have the following embedding property of the class $F_{\omega,p}(\T^d)$.
\begin{lemma}\label{lem-embedding}
	Let $0< p,q\leq \infty$ and $\omega = (\omega(k))_{k\in \Z^d}$ be a weight satisfying \eqref{eq-condition}. Then the operator $id:  F_{\omega,p}(\T^d)\hookrightarrow  F_{q}(\T^d) $ is continuous if either $p\leq q$ or $q<p$  and the series $\sum_{k\in \zz^d} \omega(k)^{-\frac{pq}{p-q}}$ converges. 
\end{lemma}
\begin{proof}
 If $q<p$ and $f\in F_{\omega,p}(\T^d)$, applying H\"older's inequality we get
\begin{align*}
\bigg(\sum_{k\in \zz^d}|\hat{f}(k)|^q\bigg)^{\frac{1}{q}} \leq \bigg(\sum_{k\in \zz^d} \omega(k)^{-\frac{pq}{p-q}}\bigg)^{\frac{p-q}{pq}} \bigg( \sum_{k\in \zz^d} |\omega(k)\hat{f}(k)|^{p}\bigg)^{{\frac{1}{p}}}.
\end{align*}
This proves the case $q<p$. The case $p\leq q$ is obvious.
\end{proof}
Our result for the best $n$-term approximation of the embedding $F_{\omega,p}(\T^d)\to F_{q}(\T^d)$ by the trigonometric system $\cT^d$ reads as follows.

%&&&&&&&&&&&&&&&&&&&&&&&&&&&&&&&&&&&&&&&&&&&&&&&&&&&&&&&&&&&&&&&&
%&&&&&&&&&&&&&&&&&&&&&&&&&&&&&&&&&&&&&&&&&&&&&&&&&&&&&&&&&&&&&&

%&&&&&&&&&&&&&&&&&&&&&&&&&&&&&&&&&&&&&&&&&&&&&&&&&&&&&&&&&&&&&&&
%&&&&&&&&&&&&&&&&&&&&&&&&&&&&&&&&&&&&&&&&&&&&&&&&&&&&&&&&&&&&&&&&

%&&&&&&&&&&&&&&&&&&&&&&&&&&&&&&&&&&&&&&&&&&&&&&&&&&&&&&&&&&&&&&&&&&&&&&&&&&&
%&&&&&&&&&&&&&&&&&&&&&&&&&&&&&&&&&&&&&&&&&&&&&&&&&&&&&&&&&&&&&&&&&&&&&&&&&&&

\begin{theorem}\label{satz1}
	Let $0<p,q\leq \infty$ and $\omega = (\omega(k))_{k\in \Z^d}$ be a weight satisfying  conditions in Lemma \ref{lem-embedding}.   Then we have
 $$
 \sigma_n\big(id: F_{\omega,p}(\T^d)\to F_{q}(\T^d),\cT^d\big)	= 	\sigma_n\big(T_\lambda:\ell_p(\N)\to \ell_q(\N),\cE\big),\quad n\in \N,
 $$
 where the  value  of $\sigma_n(T_\lambda,\cE)$ is given as in Theorem \ref{prop-sequence}.
\end{theorem}

\begin{proof} We consider the following commutative diagram
	
	\begin{equation*}
		\begin{CD}
			F_{\omega,p}(\T^d)  @ > id >> F_q(\T^d) \\
			@VV A V @AA B A\\
			\ell_p(\Z^d) @ > D_\omega >> \ell_q(\Z^d) \, ,
		\end{CD}
	\end{equation*}
	where the linear operators $A$,  $B$ and 
	$D_\omega$
	are defined as 
\begin{equation*}
	\begin{aligned}
	Af & : =  (\omega(k)\hat{f}(k))_{k\in \Z^d}\, ,  
\\
D_\omega\xi & :=   (\xi(k)/\omega(k))_{k\in \Z^d}\,  , \qquad \xi=(\xi(k))_{k\in \Z^d}
\\
(B\xi)(x)& := \sum_{k\in \Z^d} \xi_k\,  e^{{\rm i}k x}\, , \qquad x \in \T^d\, .
\end{aligned}
\end{equation*}
	It is obvious that $\|A\|=\|B\|=1$. Let $\cE^d:=\{e_k: k\in \Z^d\}$ where $e_k=(\delta_{k,l})_{l\in \Z^d}$. 
	By the property \eqref{eq-ideal} and the identity $id = B\, D_\omega \, A$ it follows
\begin{equation*}
\sigma_n\big(id: F_{\omega,p}(\T^d)\to F_{q}(\T^d),\cT^d\big)
\leq \sigma_n\big(D_\omega:\ell_p(\Z^d)\to \ell_q(\Z^d), \cE^d\big),\qquad n\in \N\, .
\end{equation*}
From the fact that
\begin{equation}\label{eq-zd-n}
\sigma_n\big(D_\omega:\ell_p(\Z^d)\to \ell_q(\Z^d), \cE^d\big) = \sigma_n\big(T_\lambda:\ell_p(\N)\to \ell_q(\N), \cE\big)
\end{equation}
we obtain the estimate from above. 
	Now we employ the same type of arguments with respect to the diagram

	\begin{equation*}
		\begin{CD}
			\ell_p(\Z^d) @ > D_\omega >> \ell_q(\Z^d)
			\\
			@VV A^{-1} V @AA B^{-1} A\\
			F_{\omega,p}(\T^d)  @ > id >> F_q(\T^d)\,.
		\end{CD}
	\end{equation*}

	It is easy to see that the operators $A$ and $B$ are invertible and  that $\|A^{-1}\|=\|B^{-1}\|=1$.
	As above we conclude 
\begin{equation*}
\sigma_n\big(D_\omega:\ell_p(\Z^d)\to \ell_q(\Z^d), \cE^d\big) \leq \sigma_n\big(id: F_{\omega,p}(\T^d)\to F_{q}(\T^d),\cT^d\big)
,\qquad n\in \N\, .
\end{equation*}
Now the estimate from below follows from \eqref{eq-zd-n}.
\end{proof}

%&&&&&&&&&&&&&&&&&&&&&&&&&&&&&&&&&&&&&&&&&&&&&&&&&&&&&&&&&&&&&&&&&&&&&&&&&&&
%&&&&&&&&&&&&&&&&&&&&&&&&&&&&&&&&&&&&&&&&&&&&&&&&&&&&&&&&&&&&&&&&&&&&&&&&&&&

We need following auxiliary results.
\begin{lemma}\label{lem:auxi-1} \begin{description}
	\item[(i)] 	Let $s>0$, $a > 1$, and $\beta\geq 0$. Then we have
	\begin{equation*}
		\lim\limits_{n\to \infty} \, \int_{\frac{a}{n}}^1 y^{s}\bigg(\frac{\ln n}{\ln(yn)}\bigg)^{\beta}\dd y =\frac{1}{s+1}.
	\end{equation*}
\item[(ii)]	Let $s>1$, $\beta \ge 0$. Then we have $$\lim_{n\to \infty} \int_1^{+\infty} \dfrac{1}{t^s} \bigg(\frac{\ln (nt)}{\ln n}\bigg)^\beta \dd t=\dfrac{1}{s-1}.$$
\end{description}	
\end{lemma}
\begin{proof}The first statement was proved in \cite{NNS20}. We prove the second one with concentration on the case $\beta>0$ since the case $\beta=0$ is obvious. We consider the sequence of functions
$$f_n(t)=\dfrac{1}{t^s}\bigg(\frac{\ln (nt)}{\ln n}\bigg)^\beta,\ \ \ t \ge 1,\ \ n\in \N.$$
Clearly, this sequence converges pointwise to $f(t)=\frac{1}{t^s}$. For $n \ge 3$, from the inequality $(x+y)^\beta \le C_\beta (x^\beta +y^\beta)$, for some $C_\beta>0$, we derive:
\begin{equation*}
    f_n(t) =\frac{1}{t^s}\bigg(1+\dfrac{\ln t}{\ln n}
    \bigg)^\beta< \frac{1}{t^s}\big(1+\ln t
    \big)^\beta \leq C_\beta \dfrac{1}{t^s}\big(1+(\ln t) ^\beta\big):=g(t).
\end{equation*}
Since $g(t)$ is integrable on $[1,+\infty)$, the desired result follows from Lebesgue's dominated convergence theorem.
\end{proof}

The asymptotic constants of best $n$-term approximation widths of embeddings of the classes $F_{\omega,p}(\T^d)$ in   $F_q(\T^d)$ are given in the following theorem.
\begin{theorem}\label{thm-general} 
Let $s>0$,  $\beta\geq 0$ and let $\omega$  be a given  weight. Assume that there exists 
$C\in \R$ such that
 \begin{equation}\label{eq-assumption}
 \lim_{n \to \infty} \, \frac{\lambda_n}{n^{-s}(\ln n)^{\beta} \, }=   \lim_{n \to \infty} \, \frac{a_n\big(id: F_{\omega,2}(\T^d)\to L_{2}(\T^d) \big)}{n^{-s}(\ln n)^{\beta}} = C\, .
 	\end{equation}
\begin{description}
	\item[(i)] If $0<p\leq q\leq \infty$ we have 
	\begin{equation*}
		\lim\limits_{n\to \infty} \frac{ \sigma_n\big(id: F_{\omega,p}(\T^d) \to F_q(\T^d),\cT^d\big)}{n^{-s-\frac{1}{p}+\frac{1}{q}}(\ln n)^{\beta}}
		= \frac{(s+\frac{1}{p}-\frac{1}{q})^{s+\frac{1}{p}-\frac{1}{q}}}{(s+\frac{1}{p})^{s}} \frac{p^{\frac{1}{p}}}{q^{\frac{1}{q}}}  C\, .
	\end{equation*}
	If $q=\infty$ and/or $p=\infty$, the asymptotic constant is understood as the limit of the right-hand side when $q\to \infty$ and/or $p\to \infty$.
\item[(ii)] If $0<q<p<\infty$ and $s>\frac{1}{q}-\frac{1}{p}$ we have 
\begin{equation*}
	\lim\limits_{n\to \infty} \frac{ \sigma_n\big(id: F_{\omega,p}(\T^d) \to F_q(\T^d),\cT^d\big)}{n^{-s-\frac{1}{p}+\frac{1}{q}}(\ln n)^{\beta}}
	=\bigg(\frac{s}{s+\frac{1}{p}} \bigg)^s \bigg(\dfrac{\frac{1}{q}}{s+\frac{1}{p}-\frac{1}{q}} \bigg)^{\frac{1}{q}-\frac{1}{p}}C .
\end{equation*}
\end{description}
\end{theorem}

 \begin{proof}  We prove the case $0<p,q<\infty$. The cases $p=\infty$ and/or $q=\infty$ are carried out similarly with slight modification.  
In this proof, for simplicity we denote
$$
\sigma_n:= \sigma_n\big(id: F_{\omega,p}(\T^d) \to F_q(\T^d),\cT^d\big).
$$
{\it Step 1.} We need some preparations. Assumption \eqref{eq-assumption} indicates that for any $\varepsilon>0$ there exists $n_1:=n_1(\varepsilon)\in \N$ such that for  $k>n_1$ we have
	\begin{equation}\label{eq-sigma-epsilon}
		\bigg|\frac{\lambda_k}{ k^{-s} (\ln k)^{\beta}}   - C\bigg| \leq \varepsilon \quad \Longleftrightarrow \quad 
		C -\varepsilon \leq \frac{\lambda_k}{k^{-s} (\ln k)^{\beta}} \leq \varepsilon + C.
	\end{equation}
Denote by $n_2=n_2(p,s)\in \N$  from which the function $\psi(t)=t^{ps}(\ln t)^{-p\beta}$ is   increasing. Then for  $m>n_0:=\max\{ n_1,n_2\}$ we have
	\begin{equation*}
		\sum_{k=1}^{m}\lambda_k^{-p} = \sum_{k=1}^{n_0}\lambda_k^{-p} + \sum_{k=n_0+1}^{m}\lambda_k^{-p}
		\leq 
		\sum_{k=1}^{n_0}\lambda_k^{-p} +  \dfrac{1}{( C -\varepsilon )^{p}}
		 \sum_{k=n_0+1}^m k^{ps} (\ln k)^{-p\beta}\,. 
	\end{equation*} 
	Estimating the  summation by an integral and afterwards changing variable $y=\frac{t}{m+1}$ we find
	\begin{equation*}
		\begin{aligned}
			\sum_{k=n_0+1}^m k^{ps} (\ln k)^{-p\beta}
			&
			\leq   \int_{n_0+1}^{m+1} t^{ps} (\ln t)^{-p\beta} \dd t
 = \frac{(m+1)^{ps+1}}{(\ln (m+1))^{p\beta}}  \int_{\frac{n_0+1}{m+1}}^1 y^{ps}\bigg(\frac{\ln(m+1)}{\ln(y(m+1))}\bigg)^{p\beta}\dd y\,.
		\end{aligned}
	\end{equation*}
By  Lemma \ref{lem:auxi-1} (i) we can choose $n_3>n_0$ such that for 
$m\geq n_3$ we have 
	\begin{equation*}
		\int_{\frac{n_0+1}{m+1}}^1 y^{ps}\bigg(\frac{\ln(m+1)}{\ln(ym+y)}\bigg)^{p\beta}\dd y \leq \frac{1+\varepsilon}{ps+1} 
		\qquad
		\text{and}
		\qquad  
		\frac{(\ln (m+1))^{p\beta}}{(m+1)^{ps+1}}  \sum_{k=1}^{n_0}\lambda_k^{-p}  \leq \varepsilon
	\end{equation*}
which leads to
	\begin{equation}\label{ws-19}
		\sum_{k=1}^{m}\lambda_k^{-p}  \leq \frac{(m+1)^{ps+1}}{(\ln (m+1))^{p\beta}} \bigg( \varepsilon +   \frac{1+\varepsilon}{( C -\varepsilon )^{p}(ps+1)}  \bigg)
	\end{equation}
	for  $m\geq n_3$. 
	Similarly we have
\begin{equation*} 
	\begin{aligned}
		\sum_{k=n_0+1}^m k^{ps} (\ln k)^{-p\beta}
		& 
		\geq \int_{n_0}^m \frac{t^{ps}}{(\ln t)^{p\beta}}  \dd t 
		\geq  \frac{m^{ps+1}}{(\ln m)^{p\beta}} \int_{\frac{n_0 }{m}}^1 y^{ps}\bigg(\frac{\ln m}{\ln(ym)}\bigg)^{p\beta}\dd y
		\geq  \frac{m^{ps+1}}{(\ln m)^{p\beta}} \frac{1-\varepsilon}{ps+1}\,,
	\end{aligned}
\end{equation*} 	
which implies
	\begin{equation}\label{ws-20}
	\sum_{k=1}^{m}\lambda_k^{-p} \geq \sum_{k=n_0+1}^{m}\lambda_k^{-p} \geq \dfrac{1}{( C+\varepsilon )^{p}}
	\sum_{k=n_0+1}^m k^{ps} (\ln k)^{-p\beta} \geq \frac{m^{ps+1}}{(\ln m)^{p\beta}}  \cdot   \frac{1-\varepsilon}{( C +\varepsilon )^{p}(ps+1)}  .
\end{equation}	
{\it Step 2.} Proof of the case $0< p\leq q<\infty$. From \eqref{ws-20} we get 
	\begin{equation*} 
		\begin{aligned}
	\frac{m-n}{\big(\sum_{k=1}^m\lambda_k^{-p}\big)^{\frac{q}{p}}} 
			& \leq
			\frac{(m-n)(\ln m)^{q\beta}}{m^{qs+\frac{q}{p}} 
			}	   \bigg(\frac{(C+\varepsilon)^{p}(ps+1)}{1-\varepsilon}\bigg)^{\frac{q}{p}}\,. 
		\end{aligned}
	\end{equation*} 
Considering the function 
$	g(t) :=\frac{t-n}{t^{qs+\frac{q}{p}}}(\ln t)^{q\beta}, \  t\in [n,\infty)\,,
$
we have
	\begin{equation*}
		g'(t)= \bigg(\frac{-t\big(qs+\frac{q}{p}-1\big)+n\big(qs+\frac{q}{p}\big)}{t^{qs+\frac{q}{p}+1}}\bigg)(\ln t)^{q\beta} + \bigg(\frac{t-n}{t^{qs+\frac{q}{p}+1}}\bigg)q\beta (\ln t)^{q\beta-1}\,.
	\end{equation*}
	We put 
\begin{equation*}
	\begin{aligned}
	f(t):= \bigg[-t\bigg(qs+\frac{q}{p}-1\bigg)+n\bigg(qs+\frac{q}{p}\bigg)\bigg]\ln t  + \big(t-n\big)q\beta  \, , \qquad t \in [n,\infty).
	\end{aligned}
\end{equation*}
Then  $g'(t)=0$ is equivalent to $f(t)=0$. We have 
	\begin{equation*} 
		\begin{aligned}
		f'(t)& = -\Big(qs+\frac{q}{p}-1\Big)(\ln t +1) +\frac{n(qs+\frac{q}{p})}{t}+q\beta 
	  <-\Big(qs+\frac{q}{p}-1\Big)\ln t  +1+q\beta\, .
		\end{aligned}
	\end{equation*}
	This implies $f'(t) < 0 $ if $t>e^{(1+q\beta)/(qs+q/p-1)}$.
	Observe that 
\begin{align*}
	f\bigg( n+\frac{n}{qs+\frac{q}{p}-1}  \bigg)      > 0,\qquad 
f\bigg( n+\frac{2n}{qs+\frac{q}{p}-1}    \bigg)       <0\quad \text{and}\quad 
\end{align*}
	for $n\geq n_4$ depending only on $p,q,s$ and $\beta$. 
Consequently the equation $f(t)=0$ (or $g'(t)=0$) has a unique solution belonging to the interval
$
		I_n:=\big[ n+\frac{n}{qs+\frac{q}{p}-1} ,\  n+\frac{2n}{qs+\frac{q}{p}-1}  \big].
$
	From this we deduce 
	\begin{equation*} 
	\begin{aligned}
	\sigma_n^q=	\sup_{m\geq n} \Bigg(	\frac{m-n}{\big(\sum_{k=1}^m\lambda_k^{-p}\big)^{\frac{q}{p}}} \Bigg) 
		&	\le   	\sup_{t\in I_n} \bigg(	\frac{(t-n)(\ln t)^{q\beta}}{t^{qs+\frac{q}{p}} 
			}\bigg)\bigg(\frac{(C+\varepsilon)^{p}(ps+1)}{1-\varepsilon}\bigg)^{\frac{q}{p}}
	\end{aligned}
\end{equation*}
which leads to
	\begin{equation*} 
		\begin{aligned}
			\frac{\sigma_n^q}{n^{1-qs-\frac{q}{p}}(\ln n)^{q\beta}}	
			&	\le   	\sup_{t\in I_n} \bigg(	\frac{(t-n)(\ln t)^{q\beta}}{n(tn^{-1})^{qs+\frac{q}{p}} 
				(\ln n)^{q\beta}}\bigg)\bigg(\frac{(C+\varepsilon)^{p}(ps+1)}{1-\varepsilon}\bigg)^{\frac{q}{p}}\,
			\\
			&	\leq  
			\sup_{ t\in \R,\, t\geq n} \bigg(	\frac{t-n }{n(tn^{-1})^{qs+\frac{q}{p}} }\bigg)\bigg(\frac{(C+\varepsilon)^{p}(ps+1)}{(1-\varepsilon)^2}\bigg)^{\frac{q}{p}}
		\end{aligned}
	\end{equation*} 
if $n$ is large enough.
	It is easy to see that the function $h(t):= \frac{t-n}{t^{qs+q/p}}$, $t \in [n,\infty)$, attains its maximum at 
	$t_0= \big(1+\frac{1}{qs+q/p-1}\big)n$. Hence, we find
	\begin{equation} \label{ws-24b}
		\begin{aligned}
	\bigg(	\frac{\sigma_n}{n^{-s-\frac{1}{p}+\frac{1}{q}}(\ln n)^{\beta}}	  \bigg)^q
			&	\leq  
			\frac{1}{(qs+\frac{q}{p}-1)\big(1+\frac{1}{qs+q/p-1}\big)^{qs+q/p}} \bigg(\frac{(C+\varepsilon)^{p}(ps+1)}{(1-\varepsilon)^2}\bigg)^{\frac{q}{p}}    
		\end{aligned}
	\end{equation} 
if $n$ is large enough. Taking the limits $n\to \infty$  and afterwards $\varepsilon \downarrow 0$ in \eqref{ws-24b} 
we obtain the upper bound. In view of Theorem \ref{prop-sequence} (i), by choosing $m\sim  \big(1+\frac{1}{qs+q/p-1}\big)n$ we also obtain the lower bound in this case.
	\\
{\it Step 3.} Proof of the case $0<q<p<\infty$. 
Firstly, we estimate $n_*$ in \eqref{eq-q<p}.  
From \eqref{eq-n-*} we have
\begin{equation*} 
(n_*-n)\lambda_{n_*}^{-p} \leq \sum_{k=1}^{n_*} \lambda_k^{-p}.
\end{equation*}
In view of 
\eqref{eq-sigma-epsilon} and  \eqref{ws-19} we get for $n\geq n_3$
\begin{align*}
    (n_*-n)(C+\varepsilon)^{-p}\frac{n_*^{ps}}{(\ln n_*)^{p \beta}}  \le \frac{(n_*+1)^{ps+1}}{(\ln (n_*+1))^{p\beta}} \bigg( \varepsilon +   \frac{1+\varepsilon}{(C -\varepsilon )^{p}(ps+1)}  \bigg)
\end{align*}
which implies 
\begin{align*} \dfrac{n_*-n}{n_*} \le \bigg(1+\frac{1}{n_*}\bigg)^{ps+1}\bigg( \varepsilon (C+\varepsilon)^p +   \frac{(1+\varepsilon)(C+\varepsilon)^p}{( C -\varepsilon )^{p}(ps+1)}  \bigg) .
\end{align*}
Therefore, for any $\epsilon>0$, exist $N_1>0$ such that for $n>N_1$ we have 
\begin{align}\label{n*1}  \dfrac{n_*-n}{n_*} \le \dfrac{1}{ps+1}+\epsilon  \qquad\text{or} \qquad
    n_* \le \dfrac{n}{\frac{ps}{ps+1}-\epsilon}.
\end{align}

Using \eqref{eq-sigma-epsilon} and \eqref{ws-20}   the condition  $(m-n)\lambda_m^{-p} \leq \sum_{k=1}^m \lambda_k^{-p}$
  is   satisfied if $$ \frac{m^{ps+1}}{(\ln m)^{p\beta}}  \cdot   \frac{1-\varepsilon}{( C +\varepsilon )^{p}(ps+1)} \ge \dfrac{(m-n)(C-\varepsilon)^{-p}m^{ps}}{(\ln m)^{p \beta}}$$
which is equivalent to 
$$ \dfrac{m-n}{m} \le \dfrac{(1-\varepsilon)(C-\varepsilon)^p}{(C+\varepsilon)^p(ps+1)}. $$
Hence, for any $\epsilon>0$, there exists $N_2\in \N$ such that for $n>N_2$ we have $$\dfrac{(1-\varepsilon)(C-\varepsilon)^p}{(C+\varepsilon)^p(ps+1)} \ge \dfrac{1}{ps+1}-\epsilon.
$$
Therefore, the condition $
		(m-n)\lambda_m^{-p} \leq \sum_{k=1}^m \lambda_k^{-p}$ is satisfied if $$\dfrac{m-n}{m} \le \dfrac{1}{ps+1}-\epsilon\qquad \text{or}\qquad m \le \dfrac{n}{\frac{ps}{ps+1}+\epsilon}.$$
This leads to $
		    n_* \ge \frac{n}{\frac{ps}{ps+1}+\epsilon}.
$		From this and \eqref{n*1}  we deduce \begin{equation*}
		    n_* \sim \bigg(1+\frac{1}{ps}\bigg)n,\ \ \ n \to +\infty\,.
		\end{equation*}
Denoting $\alpha=\frac{pq}{p-q}$,  from \eqref{eq-q<p} we have
\begin{align} \label{esd}
\sigma_n^{\alpha}=\Bigg(\dfrac{(n_*-n)^{1/q}}{\big(\sum_{k=1}^{n_*}\lambda_k^{-p}\big)^{1/p}}\Bigg)^\alpha + \sum_{k=n_*+1}^\infty \lambda_k^{\alpha}.
\end{align}
		Using \eqref{ws-19} and \eqref{ws-20} again we get \begin{equation*}
		    \sum_{k=1}^{n_*}\lambda_k^{-p} \sim \dfrac{1}{C^p(ps+1)}\dfrac{n_*^{ps+1}}{(\ln n_*)^{p \beta}} \sim \dfrac{1}{C^p(ps+1)} \dfrac{(1+\frac{1}{ps})^{ps+1}n^{ps+1}}{(\ln n)^{p \beta}}, \ \ n \to +\infty.
		\end{equation*}
		Therefore, the first term in \eqref{esd} can be estimated: 
\begin{equation}\label{eq-first}
	\begin{aligned}
	\Bigg(	    \dfrac{(n_*-n)^{\frac{1}{q}}}{\big(\sum_{k=1}^{n_*}\lambda_k^{-p}\big)^{\frac{1}{p}}} \Bigg)^\alpha 
	&\overset{  n \to +\infty}{\sim} \Bigg(  \dfrac{C(ps+1)^{\frac{1}{p}}(\ln n)^{\beta}}{(1+\frac{1}{ps})^{s+\frac{1}{p}}  n^{s+\frac{1}{p}}} \Big(\frac{n}{ps}\Big)^{\frac{1}{q}} \Bigg)^\alpha   =\frac{C^\alpha}{ps} \bigg(\frac{ps}{ 1+ps}\bigg)^{s\alpha} \frac{(\ln n)^{\alpha\beta}}{n^{\alpha(s+\frac{1}{p}-\frac{1}{q})}}  .
	\end{aligned}
\end{equation}
Now, we estimate the second term in \eqref{esd}. Observe that   $f(t)=t^{-s\alpha}(\ln t)^{\alpha \beta}$ is a decreasing function when $t\geq t_0$ for some $t_0>0$. Hence,  when $n$ is large enough, in view of \eqref{eq-sigma-epsilon} we can bound 
		\begin{align}\label{s1} \sum_{k=n_*+1}^\infty \lambda_k^{\alpha} &\le (C+\varepsilon)^{\alpha}\sum_{k=n_*+1}^\infty \frac{(\ln k)^{\alpha \beta}}{k^{s \alpha}}  \le (C+\varepsilon)^\alpha \int_{n_*}^{+\infty} \frac{(\ln t)^{\alpha \beta}}{t^{s \alpha}}  \dd t \notag 
			\\
		&=(C+\varepsilon)^\alpha \frac{(\ln n_*)^{\alpha \beta}}{n_*^{s \alpha -1}}   \int_1^{+\infty} \dfrac{1}{t^{s \alpha}} \bigg(\dfrac{\ln n_* t}{\ln n_*}\bigg)^{\alpha \beta } \dd t.
		\end{align}
Similarly, we also have the estimate
		\begin{align} \label{s2}
		  \sum_{k=n_*+1}^\infty \lambda_k^{\alpha} &\ge (C-\varepsilon)^\alpha \sum_{k=n_*+1}^\infty \frac{(\ln k)^{\alpha \beta}}{k^{s \alpha}} \ge  (C-\varepsilon)^\alpha \int_{n_*+1}^{+\infty}\frac{ (\ln t)^{\alpha \beta}}{t^{s \alpha}} \dd t \notag \\
		  &=(C-\varepsilon)^\alpha \frac{(\ln (n_*+1))^{\alpha \beta}}{(n_*+1)^{s \alpha-1}}  \int_1^{+\infty} \dfrac{1}{t^{s \alpha}}  \bigg(\dfrac{\ln (n_*+1) t}{\ln (n_*+1)}\bigg)^{\alpha \beta } \dd t.
		\end{align}
Note that the condition $s>\frac{1}{q}-\frac{1}{p}$ implies $s\alpha>1$.	Using Lemma \ref{lem:auxi-1} (ii), from \eqref{s1} and \eqref{s2} we get
		\begin{align*}
		 \sum_{k=n_*+1}^\infty \lambda_k^{\alpha} 
		 & \overset{  n \to +\infty}{\sim}  \dfrac{1}{s \alpha -1}C^\alpha   \bigg(1+\frac{1}{ps}\bigg)^{1-s \alpha}\frac{(\ln n)^{\alpha \beta} }{ n^{s \alpha-1}}
		 \\
		&\ \ \   =  \dfrac{p-q}{spq-p+q}\bigg(1+\frac{1}{ps}\bigg)C^\alpha    \bigg( \frac{ps}{ps+1}\bigg)^{\alpha s} \frac{(\ln n)^{\alpha \beta} }{ n^{ \alpha(s+\frac{1}{p}-\frac{1}{q})}} .
		\end{align*}
From this and \eqref{eq-first} we finally obtain 
		\begin{align*}
		    \sigma_n^\alpha & \overset{  n \to +\infty}{\sim}
		    \bigg[\frac{1}{ps} + \dfrac{p-q}{spq-p+q}\Big(1+\frac{1}{ps}\Big)\bigg]  C^\alpha \bigg(\frac{ps}{ps+1}\bigg)^{s \alpha}   \frac{(\ln n)^{\alpha \beta} }{ n^{ \alpha(s+\frac{1}{p}-\frac{1}{q})}}
		    \\
		    &\ \ \  = \frac{p}{spq-p+q} C^\alpha \bigg(\frac{ps}{ps+1}\bigg)^{s \alpha}    \frac{(\ln n)^{\alpha \beta} }{ n^{ \alpha(s+\frac{1}{p}-\frac{1}{q})}}
		\end{align*}
which proves the second statement.
		\end{proof}

%&&&&&&&&&&&&&&&&&&&&&&&&&&&&&&&&&&&&&&&&&&&&&&&&&&&&&&&&&&&&&&&&&&&&&&&&&&&&&&&&&&&&&&&&&&&&&&&
%&&&&&&&&&&&&&&&&&&&&&&&&&&&&&&&&&&&&&&&&&&&&&&&&&&&&&&&&&&&&&&&&&&&&&&&&&&&&&&&&&&&&&&&&&&&&&&

%&&&&&&&&&&&&&&&&&&&&&&&&&&&&&&&&&&&&&&&&&&&&&&&&&&&&&&&&&&&&&&&&&&&&&&&&&&&
%&&&&&&&&&&&&&&&&&&&&&&&&&&&&&&&&&&&&&&&&&&&&&&&&&&&&&&&&&&&&&&&&&&&&&&&&&&&

\section{Best $n$-term approximation of function spaces with mixed smoothness}
\label{sec-mixed}

%&&&&&&&&&&&&&&&&&&&&&&&&&&&&&&&&&&&&&&&&&&&&&&&&&&&&&&&&&&&&&&&&&&&&&&&&&&&
%&&&&&&&&&&&&&&&&&&&&&&&&&&&&&&&&&&&&&&&&&&&&&&&&&&&&&&&&&&&&&&&&&&&&&&&&&&&

In this section we shall apply the  result in Section \ref{sec-F-Class} to the  family of weights 
\begin{align*} 
\omega_{s,r}(k) := & \prod_{i=1}^d\big(1+|k_i|^r\big)^{s/r}\, , \qquad 0 < r < \infty\, , \\
\omega_{s,r}(k)  := & \prod_{i=1}^d \max (1,|k_i|)^{s}\, , \qquad  r = \infty\, , 
\end{align*}
$k \in \Z^d$, where the parameter $s$ satisfies  $0<s<\infty$. 
We shall use the notation   ${H}_{\mix}^{s,r}(\T^d):={F}_{\omega_{s,r},2}(\T^d)$ and $\mathcal{A}_{\mix}^{s,r}(\T^d):=F_{\omega_{s,r},1}(\T^d)$, respectively. The classes ${H}_{\mix}^{s,r}(\T^d)$ 
 are called periodic Sobolev spaces with mixed smoothness and well-known in approximation theory, see, e.g., \cite{NoWo08,NoWo10,DTU18B}. The classes $\mathcal{A}_{\mix}^{s,r}(\T^d)$ are the weighted Wiener algebras. These spaces have been studied extensively recently in \cite{KPV15,BKUV17,KV19,NNS20}.
In both spaces ${H}_{\mix}^{s,r}(\T^d)$ and $\mathcal{A}_{\mix}^{s,r}(\T^d)$, for different $r$, we obtain the same sets of functions.
A change of the parameter $r$ leads to a change of the quasinorm only.

Let $m\in \N$. We define  the space $H^m_{\mix}(\T^d)$ to be the collection of all functions 
$f \in L_2 (\T^d)$ such that all distributional derivatives $D^\alpha f$ with $\alpha=(\alpha_1,\ldots,\alpha_d)$, $\max_{j=1,\ldots,d}\alpha_j\leq m$ 
belong to $L_2 (\T^d)$. The space $H^m_{\mix}(\T^d)$ is equipped with the norm
\begin{equation*}
\big\|\, f\, \big\|_{H^m_{\mix}(\T^d)}  := \Bigg(\sum_{\alpha =(\alpha_1,\ldots,\alpha_d)\in \N_0^d\atop \alpha_j \leq m,\, j=1,\ldots,d} \, \big\|\, D^{\alpha}f\big\|_{L_2(\T^d)}^2\Bigg)^{1/2}\, .
\end{equation*}  
Then $H^m_{\mix}(\T^d) = H^{m,r}_{\mix}(\T^d)$ for all $r$ in the sense of equivalent quasinorms.
If $m=1$, then we have 
$
\| \cdot  \|_{H^{1,2}_{\mix}(\T^d)}= \| \cdot\|_{H^1_{\mix}(\T^d)}.
$
If $m \ge 2$, then the norm $\| \cdot\|_{H^m_{\mix}(\T^d)}$ itself does not belong to the family of norms
$\| \cdot\|_{H^{m,r}_{\mix}(\T^d)}$, $0 < r \le \infty$.
But the choice $r=2m$ leads to the following standard norm  
\begin{equation*}
\big\|\, f\, \big\|_{H^{m,2m}_{\mix}(\T^d)}= \Bigg(\sum_{\alpha \in \{0,m\}^d} \big\| \, D^{\alpha}f\, \big\|_{L_2
	(\T^d)}^2\Bigg)^{1/2} \, ,
\end{equation*}
see \cite{KSU15}.

%But in our final Section \ref{ex3} we will study the embedding ${H}^{s,2}_{\mix} (\T^d) \hookrightarrow H^{1,2}(\T^d)$.

%&&&&&&&&&&&&&&&&&&&&&&&&&&&&&&&&&&&&&&&&&&&&&&&&&&&&&&&&&&&&&&&&&&&&&&&&&&&&&&&&&&
%&&&&&&&&&&&&&&&&&&&&&&&&&&&&&&&&&&&&&&&&&&&&&&&&&&&&&&&&&&&&&&&&&&&&&&&&&&&&&&&&&&&

%&&&&&&&&&&&&&&&&&&&&&&&&&&&&&&&&&&&&&&&&&&&&&&&&&&&&&&&&&&&&&&&&&&&&&&&&&&&&&&&&&&&&
%&&&&&&&&&&&&&&&&&&&&&&&&&&&&&&&&&&&&&&&&&&&&&&&&&&&&&&&&&&&&&&&&&&&&&&&&&&&&&&&&&&&&&

%&&&&&&&&&&&&&&&&&&&&&&&&&&&&&&&&&&&&&&&&&&&&&&&&&&&&&&&&&&&&&&&&&&&&&&&&&&&&&&&&&&
%&&&&&&&&&&&&&&&&&&&&&&&&&&&&&&&&&&&&&&&&&&&&&&&&&&&&&&&&&&&&&&&&&&&&&&&&&&&&&&&&&&&

%&&&&&&&&&&&&&&&&&&&&&&&&&&&&&&&&&&&&&&&&&&&&&&&&&&&&&&&&&&&&&&&&&&&&&&&&&&&&&&&&&&&&
%&&&&&&&&&&&&&&&&&&&&&&&&&&&&&&&&&&&&&&&&&&&&&&&&&&&&&&&&&&&&&&&&&&&&&&&&&&&&&&&&&&&&&

Let
$\lambda=(\lambda_n)_{n\in \N}$ denote the  non-increasing rearrangement of the sequence  $(1/\omega_{s,r}(k))_{k\in \Z^d}$. That leads to $\lambda_n=a_n\big(id: H^{s,r}_{\mix}(\T^d)\to L_2(\T^d)\big)$.
We recall a  result obtained in \cite{KSU15}.
\begin{prop}\label{alt}
Let $0<s<\infty$ and $0<r\leq \infty$. 
Then it holds
\begin{equation*} 
\lim\limits_{n\to \infty} \frac{\lambda_n}{n^{-s}(\ln n)^{s(d-1)}} = \lim\limits_{n\to \infty} \frac{a_n\big(id: H^{s,r}_{\mix}(\T^d)\to L_2(\T^d)\big)}{n^{-s}(\ln n)^{s(d-1)}}= 
\bigg( \frac{2^d}{(d-1)!}\bigg)^s\,.
\end{equation*}
\end{prop}

From this and Theorem \ref{thm-general} we get the following.

 \begin{theorem}\label{thm:a-mix} Let $0<s<\infty$ and $0<r\leq \infty$.  
 	Then it holds
 	\begin{equation*} 
 		\lim\limits_{n\to \infty} \frac{\sigma_n\big(id: H^{s,r}_{\mix}(\T^d)\to L_2(\T^d),\cT^d\big)}{n^{-s}(\ln n)^{s(d-1)}}=\frac{s^s}{(s+\frac{1}{2})^s} \bigg( \frac{2^d}{(d-1)!}\bigg)^s\,
 	\end{equation*}
and if $s>1/2$
\begin{equation*} 
	\lim\limits_{n\to \infty} \frac{\sigma_n\big(id: H^{s,r}_{\mix}(\T^d)\to \mathcal{A}(\T^d),\cT^d\big)}{n^{-s+\frac{1}{2}}(\ln n)^{s(d-1)}}=\bigg(\frac{s}{s+\frac{1}{2} }\bigg)^s\bigg(\frac{1}{ s-\frac{1}{2}}\bigg)^{\frac{1}{2}}\bigg( \frac{2^d}{(d-1)!}\bigg)^s\,.
\end{equation*} 
\end{theorem}
  
 The asymptotic constants for embeddings of best $n$-term approximation widths of embedding of weighted Wiener classes  $\mathcal{A}^{s,r}_{\mix}(\T^d)$ are also obtained from Theorem \ref{thm-general}.    
 \begin{theorem} \label{thm:a-mix-2}
 	 Let $0<s<\infty$ and $0<r\leq \infty$.  
 	Then it holds
	\begin{equation*} 
		\lim\limits_{n\to \infty} \frac{ \sigma_n\big(id: \mathcal{A}_{\mix}^{s,r}(\T^d) \to L_2(\T^d),\cT^d\big)}{n^{-s-\frac{1}{2}}(\ln n)^{s(d-1)}}
		= \frac{(s+\frac{1}{2})^{s+\frac{1}{2}}}{\sqrt{2}(s+1)^s}\bigg( \frac{2^d}{(d-1)!}\bigg)^s
	\end{equation*}
and
\begin{equation*} 
		\lim\limits_{n\to \infty} \frac{ \sigma_n\big(id: \mathcal{A}_{\mix}^{s,r}(\T^d) \to \mathcal{A}(\T^d),\cT^d\big)}{n^{-s}(\ln n)^{s(d-1)}}
		= \frac{s^s}{(s+1)^s}\bigg( \frac{2^d}{(d-1)!}\bigg)^s.
	\end{equation*}
\end{theorem}

\begin{remark}Let us compare the asymptotic decay of $a_n$ and $\sigma_n$. The equivalence
$$
a_n\big(id: H^{s,r}_{\mix}(\T^d)\to L_2(\T^d)\big)\asymp\sigma_n\big(id: H^{s,r}_{\mix}(\T^d)\to L_2(\T^d),\cT^d\big)
$$
has been known with a long history, see, e.g., \cite[Chapters 4 and 7]{DTU18B} for comments.
From Theorem \ref{thm:a-mix} and \cite{NNS20} we also have
$$
\sigma_n\big(id: H^{s,r}_{\mix}(\T^d)\to \mathcal{A}(\T^d),\cT^d\big) \asymp \sigma_n\big(id: H^{s,r}_{\mix}(\T^d)\to \mathcal{A}(\T^d),\cT^d\big).
$$
However, by Theorem \ref{thm:a-mix-2} and \cite{NNS20} we find
$$
a_n\big(id: \mathcal{A}_{\mix}^{s,r}(\T^d) \to L_2(\T^d)\big) \asymp n^{\frac{1}{2}} \sigma_n\big(id: \mathcal{A}_{\mix}^{s,r}(\T^d) \to L_2(\T^d),\cT^d\big).
$$
This indicates that approximating functions in the class $\mathcal{A}_{\mix}^{s,r}(\T^d)$ by $n$-term improves the convergence rate $1/2$ compared to  linear method.
\end{remark}

%&&&&&&&&&&&&&&&&&&&&&&&&&&&&&&&&&&&&&&&&&&&&&&&&&&&&&&&&&&&&&&&&&&&&&&&&&&&&&&&&&&
%&&&&&&&&&&&&&&&&&&&&&&&&&&&&&&&&&&&&&&&&&&&&&&&&&&&&&&&&&&&&&&&&&&&&&&&&&&&&&&&&&&&

%&&&&&&&&&&&&&&&&&&&&&&&&&&&&&&&&&&&&&&&&&&&&&&&&&&&&&&&&&&&&&&&&&&&&&&&&&&&
%&&&&&&&&&&&&&&&&&&&&&&&&&&&&&&&&&&&&&&&&&&&&&&&&&&&&&&&&&&&&&&&&&&&&&&&&&&&

%&&&&&&&&&&&&&&&&&&&&&&&&&&&&&&&&&&&&&&&&&&&&&&&&&&&&&&&&&&&&&&&&&&&&&&&&&&&&&&&&&&&&&&

%&&&&&&&&&&&&&&&&&&&&&&&&&&&&&&&&&&&&&&&&&&&&&&&&&&&&&&&&&&&&&&&&&&&&&&&&&&&&&&&&&&&&&&

We are also interested in asymptotic constants of best $n$-term approximation of embeddings of function spaces with mixed smoothness into $ H^1(\T^d)$. 
Here   $H^1(\T^d)$ is equipped with the norm 
\begin{align*}
 \|\, f\, \|_{H^1(\T^d)}:&= \Bigg(\sum_{k \in \Z^d} \bigg(1+\sum_{j=1}^d |k_j|^2\bigg) |\hat{f}(k)|^2\Bigg)^{1/2} 
 = 
\Bigg(\|\, f\, \|_{L_2(\T^d)}^2 + \sum_{j=1}^d \Big\|\, \frac{\partial f}{\partial x_j}\,   \Big\|_{L_2(\T^d)}^2\Bigg)^{1/2}
\, .
\end{align*}
I.e., $H^1(\T^d)$ is the standard isotropic periodic Sobolev space with smoothness $1$. 
We define a weight $\tilde{\omega}$ by
\begin{equation}\label{eq-omega}
\tilde{\omega}(k):= \frac{\prod_{j=1}^d(1+|k_j|^2)^{s/2}}{\big(1+\sum_{j=1}^d |k_j|^2\big)^{1/2}},\qquad k=(k_1,\ldots,k_d)\in \Z^d\, .
\end{equation}
Rearranging non-increasingly the sequence $(1/\tilde{\omega}(k))_{k\in \Z^d}$ with the outcome denoted by $(\tilde{\lambda}_{n})_{n\in \N}$, 
we obtain $\tilde{\lambda}_n = a_n\big(id: H^{s,2}_{\mix}(\T^d)\to H^1(\T^d)\big)$. The asymptotic constant of $a_n\big(id: H^{s,2}_{\mix}(\T^d)\to H^1(\T^d)\big)$ was obtained recently in \cite{NNS20}.

\begin{prop}\label{thm-dk} Let $d\in \N$, $s >1$ and 
\begin{equation}\label{eq-S}
S:=\sum_{k=1}^{+\infty}\dfrac{1}{(k^2+1)^{\frac{s}{2(s-1)}}}.
\end{equation}
Then we have
$$\lim _{n \to +\infty} \frac{\tilde{\lambda}_n }{n^{1-s}}=\lim _{n \to +\infty} \frac{a_n\big(id: H^{s,2}_{\mix}(\T^d)\to H^1(\T^d)\big)}{n^{1-s}}=(2d)^{s-1}(2S+1)^{(s-1)(d-1)}.$$ 
\end{prop}

From this and Theorem \ref{thm-general} we obtain the following.

\begin{theorem}\label{thm:a-mix-h1} Let $d\in \N$, $s>1$ and $S$ be given in \eqref{eq-S}. Then it holds
	\begin{equation*}
		\lim\limits_{n\to \infty} \frac{ \sigma_n\big(id: H_{\mix}^{s,2}(\T^d) \to H^1(\T^d),\cT^d\big)}{n^{-s+1}}= \bigg(\frac{s-1}{s-\frac{1}{2}}\bigg)^{s-1} (2d)^{s-1}(2S+1)^{(s-1)(d-1)} 
	\end{equation*}
and
	\begin{equation*} 
		\lim\limits_{n\to \infty} \frac{ \sigma_n\big(id: \mathcal{A}_{\mix}^{s,2}(\T^d) \to H^1(\T^d),\cT^d\big)}{n^{-s+\frac{1}{2}}}
		= \frac{(s-\frac{1}{2})^{s-\frac{1}{2}}}{\sqrt{2}\,s^{s-1}} (2d)^{s-1}(2S+1)^{(s-1)(d-1)}.
	\end{equation*}
\end{theorem}
\begin{proof}Let $\tilde{\omega}$ be given in \eqref{eq-omega}.  We will show that
	\begin{equation}\label{eq-Amix-Fomega}
		\sigma_n\big(id: H_{\mix}^{s,2}(\T^d) \to H^1(\T^d),\cT^d \big) = \sigma_n\big(id: F_{{\tilde{\omega}},2}(\T^d)\to F_2(\T^d),\cT^d\big)
	\end{equation}
and
\begin{equation}\label{eq-Amix-Aomega}
\sigma_n\big(id: \mathcal{A}_{\mix}^{s,2}(\T^d) \to H^1(\T^d),\cT^d \big) = \sigma_n\big(id: F_{{\tilde{\omega}},1}(\T^d)\to F_2(\T^d),\cT^d\big),
\end{equation}
by using standard lifting arguments.  We consider the diagram
	\begin{equation*}
	\begin{CD}
		H_{\mix}^{s,2}(\T^d)  @ > id >> H^1(\T^d) \\
		@VV A V @AA B A\\
		F_{{\tilde{\omega}},2}(\T^d) @ > id >> F_2(\T^d) \, 
	\end{CD}
\end{equation*}
where the linear operators $A$ and $B$ 	are defined for $f\in 	H^{s,2}_{\mix}(\T^d) $ and $g\in F_2(\T^d)$ respectively by
\begin{equation*}
\widehat{Af}(k)  : =  \bigg(1+\sum_{j=1}^d |k_j|^2\bigg)^{1/2}\hat{f}(k),  
\quad 
\widehat{Bg}(k) :=  \bigg(1+\sum_{j=1}^d |k_j|^2\bigg)^{-1/2}\hat{g}(k),\ \ k=(k_1,\ldots,k_d) \in \Z^d.
\end{equation*} 
It is obvious that $\|A\|=\|B\|=1$. Now by the   property  \eqref{eq-ideal}, we obtain
$$
\sigma_n\big(id: H_{\mix}^{s,2}(\T^d) \to H^1(\T^d),\cT^d\big) \leq  \sigma_n\big(id: F_{{\tilde{\omega}},2}(\T^d)\to F_2(\T^d),\cT^d\big).
$$
The reverse inequality  follows from the modified diagram
	\begin{equation*}
	\begin{CD}
		H_{\mix}^{s,2}(\T^d)  @ > id >> H^1(\T^d)\\
		@AA {A}^{-1} A @VV {B}^{-1} V\\
	F_{{\tilde{\omega}},2}(\T^d) @ > id >> F_2(\T^d) \, .
	\end{CD}
\end{equation*}
Hence  \eqref{eq-Amix-Fomega} is proved. Proof of \eqref{eq-Amix-Aomega} is carried out similarly.
Now the assertion follows from Proposition \ref{thm-dk} and Theorem \ref{thm-general}.
\end{proof}

%&&&&&&&&&&&&&&&&&&&&&&&&&&&&&&&&&&&&&&&&&&&&&&&&&&&&&&&&&&&&&&&
%&&&&&&&&&&&&&&&&&&&&&&&&&&&&&&&&&&&&&&&&&&&&&&&&&&&&&&&&&&&&&&&&
\bibliographystyle{abbrv}

\bibliography{AllBib}
\end{document}